\documentclass[11pt]{article}
  
\usepackage{graphicx}

\usepackage{amsmath,amsfonts,amssymb}

\usepackage{color}

\newcommand{\ignore}[1]{}
\newcommand{\eqn}[1]{(\ref{#1})}
\newcommand{\If}{\mathcal{I}_f}
\newcommand{\half}{\frac{1}{2}}
\newcommand{\dt}{\Delta t}
\newenvironment{mat}{\left[ \begin{array}{cccc}}{\end{array}\right]}
\newcommand\bcm{\begin{mat}}
\newcommand\ecm{\end{mat}}
\newcommand{\vu}{\vec{u}}

\newcommand{\vz}{\vec{z}}

\newcommand{\vpsi}{\vec{\psi}}
\newcommand{\grad}{\nabla}

\numberwithin{equation}{section}

\begin{document}

\title{Towards Adaptive Simulations of Dispersive Tsunami 
Propagation from an Asteroid Impact}

\author{Marsha J. Berger\thanks{Courant Institute, New York University, 251 Mercer St., NY,
NY 10012, U.S.A. Also Flatiron Institute, 162 5th Ave., NY, NY 10010, \texttt{berger@cims.nyu.edu}.}
\and
Randall J. LeVeque\thanks{Department of Applied Math, University of Washington,
Seattle, WA 98195, U.S.A. Also HyperNumerics LLC, Seattle, Washington, 
\texttt{rjl@uw.edu}.}}

\date{September 8, 2021}

\maketitle

\begin{center}
{\large  Submitted to the Proceedings of the \\
International Congress of Mathematicians, 2022}
\end{center}

\vskip 5pt

{\bf Abstract.}
The long-term goal of this work is the development of high-fidelity simulation 
tools for dispersive tsunami propagation. A dispersive 
model is especially important for short wavelength phenomena such as an asteroid 
impact into the ocean, and is also important
in modeling other events where the simpler shallow water equations are insufficient.
Adaptive simulations are crucial to bridge the scales from deep ocean
to inundation, but have difficulties with the implicit system of equations that
results from dispersive models.  We propose a fractional step scheme that advances
the solution on separate patches with different spatial resolutions and time steps.
We show a simulation with 7 levels of adaptive meshes and 
onshore inundation resulting from a simulated asteroid impact off 
the coast of Washington.
Finally, we discuss a number of open research questions that need to be resolved 
for high quality simulations.

\section{Introduction}\label{sec:intro}

Many steps are required in modeling a tsunami arising from an asteroid 
impact in the ocean. The impact itself forms a crater that drives the 
eventual tsunami creation. Modeling this requires a complex three-dimensional 
multi-physics hydrocode, since there are many physical processes and time 
scales. Once the tsunami has formed, it propagates hundreds or thousands of 
kilometers across the ocean. When the shoreline is reached, the ultimate 
goal is modeling the inundation risk to coastal populations and important 
infrastructure at a much smaller spatial scale (typically 10 meters or less).

This work addresses the last two steps, the long-distance propagation and 
coastal inundation. The goal is a high-fidelity model that can accurately 
determine the inundation risk for particular sites using available 
bathymetric data sets. Since large scale ocean simulations are so 
compute-intensive, for many tsunami modeling problems the two-dimensional 
depth-averaged Shallow Water Equations (SWE) are used for the
propagation step.  These equations assume the wavelength is long relative
to the depth of the ocean, as is typical for tsunamis generated by large
earthquakes.   However, these non-dispersive equations are often insufficient for 
short-wavelength asteroid-generated tsunamis, giving inaccurate results for 
both tsunami travel time and maximum shoreline run-in, as noted in our
own work \cite{NASA-AGT-2018} and in several other studies, e.g.
\cite{KorycanskyLynett2007,RobertsonGisler2019,WardAsphaug:2000}.
This is also the case for landslide-generated tsunamis and other short wavelength
phenomena; see e.g. \cite{Glimsdal2013}.
Shorter waves experience significant dispersion (waves with different
periods propagate with different speeds), while the hyperbolic SWE are 
non-dispersive.  Dispersive depth-averaged equations can be obtained
by retaining more terms when reducing from the three-dimensional Euler
equations to two space dimensions, giving some form of "Boussinesq equations".
The additional terms involve higher-order derivatives (typically third
order), and several different models have been proposed. 
When solved numerically, these equations generally require implicit methods
in order to remain stable with physically reasonable time steps.  By
contrast, the hyperbolic SWE involve only first-order derivatives and 
explicit methods are commonly used.  

Our numerical model is based on the
GeoClaw softwave (part of the open source Clawpack software project
\cite{clawpack}), which has been heavily used and well-validated for modeling
earthquake-generated tsunamis using the SWE 
\cite{BergerGeorgeEtAl2011a,geoclaw-nthmp-results:2011,LeVequeGeorgeEtAl2011}.
The numerical methods used are high-resolution, shock-capturing finite volume
methods based on Riemann solvers, a standard approach for nonlinear 
hyperbolic problems \cite{fvmhp}.  In the case of GeoClaw, additional
features are included to make the methods ``well balanced'' so that the
steady state of an ocean at rest is preserved.  Moreover, the shoreline
is represented as an interface between wet cells and dry cells and robust
Riemann solvers allow the determination of the fluxes at these interfaces.
The wet/dry status of a cell can change dynamically
as the tsunami advances onshore or retreats.
This software implements adaptive mesh refinement (AMR), critical for solving
problems with vastly different spatial scales from trans-ocean propagation
to community-level inundation modeling.  However, the patch-based
AMR algorithms are based on the use of explicit solvers, and the extension
of the GeoClaw software to also work with implicit solvers for Boussinesq equations
has been a major part of this project.  The basic approach used can also be
used more generally with the AMR version of Clawpack, which has many potential
applications to other wave propagation problems when dispersive or dissipative
(e.g. second-order derivative) terms are included and implicit solvers are needed.

\section{Overview of Approach }\label{sec:overview}

We are building on the work of \cite{Kim2017}, in which the 
Boussinesq equations described in the next section were solved using an extension
of GeoClaw, but only for the case of
a single grid resolution, without the AMR capability.
These equations have the form of the two-dimensional SWE with the addition
of ``source terms'' involving third order derivatives.
Equations of this type can often be solved by fractional step or splitting
methods: first advance the solution by solving the hyperbolic
shallow water equations, and then advance the solution using
terms associated with the higher order 
derivatives (or in the opposite order, as we have found to be advantageous). 
For the SWE part, we can use the standard GeoClaw solver.  The 
third derivative terms require an implicit step as described below.
We are currently using a sparse linear system solver
called Pardiso \cite{pardiso-7.2b}; 
several other groups have found that multigrid 
works nicely as well.

\begin{figure}[t]
\centerline{\includegraphics[height=2.1in]{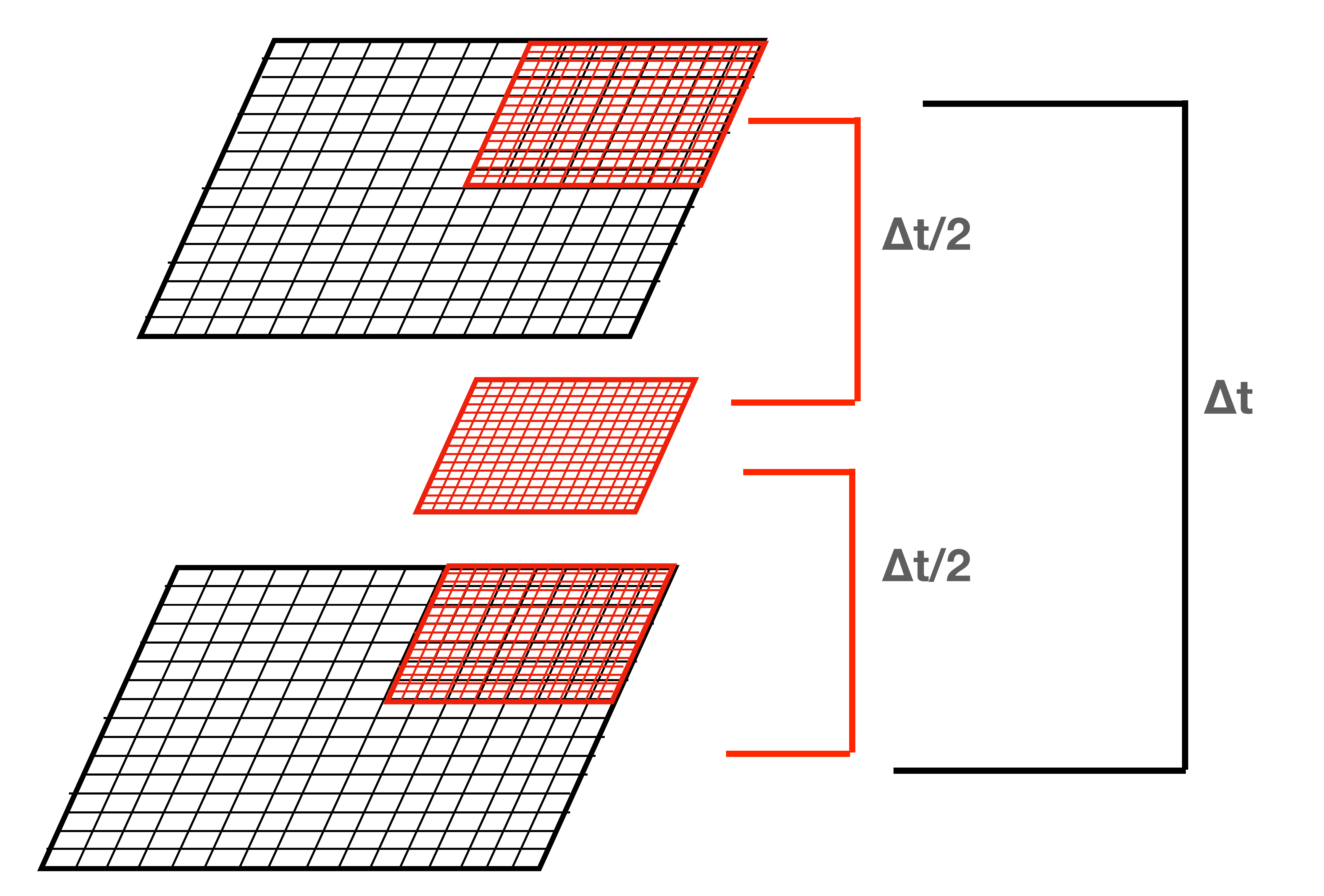}}
\caption{Figure shows a coarse grid with coarse cells outlined in the base
grid, and a fine patch refined by a factor of two, also with cells
outlined. The fine grid time step is half the coarse grid step. Before
the fine grid takes a step, ghost cells are needed to complete the
stencil.\label{fig:amrfig}}
\end{figure}

The difficulty in the solution algorithm comes from the combination of 
adaptive mesh refinement and implicit solvers.  Adaptive mesh refinement 
is critical in bridging the scales of oceanic tsunamic propagation, which 
typically needs resolution on the order of kilometers, and inundation 
modeling, where a resolution on the order of 10 meters or less is required.  
The patch-based mesh refinement in GeoClaw refines in time as well as space, 
in order to satisfy the time step stability constraint of explicit methods.
If a patch is refined in space by a factor of 4, then we also typically
refine in time by the same factor (as required by the CFL condition for
explicit methods) and so 4 time steps are taken for the 
fine patch to ``catch up'' to the coarse patch in time. The fine patch thus needs 
to interpolate ghost cells that fill out the stencil at each intermediate time step. 
Figure \ref{fig:amrfig} indicates this schematically for refinement by 2.

In our approach, the solution of the implicit equations is stored as additional 
elements in the solution vector, increasing the number of equations in two 
horizontal dimensions from 3 to 5. 
We also reverse the typical splitting order and perform the implicit 
solve first.  At the initial and final times during this time step, the ghost cells 
are still interpolated in space from the coarse grid but do not need interpolation 
in time. The ghost cell values at intermediate times on the fine grid
are interpolated from the coarse grid values at times $n$ and $n+1$.

There are other issues to consider when combining Boussinesq and SWE in a 
single solver.  The Boussinesq equations give a better model of dispersive
waves over some regime, but lack any wave breaking mechanism. As large waves
approach shore this can lead to very large magnitude solitary waves that should
break.  The nonlinear SWE perform better at this point; a shock
wave develops that is a better representation of the turbulent bore formed by
a breaking wave.  Very large waves, such as those that might be formed by an
asteroid impact, can undergo shoaling far out on the continental shelf and
dissipate some of their energy in this manner (the van Dorn effect 
\cite{KorycanskyLynett}).
(Shoaling refers to the modification of wave heights when the wave enters
shallower water, when the wave steepens and becomes higher frequency).
So it is important to transition from using the Boussinesq equations in deep
water to SWE closer to shore, based on some breaking criterion.
Also the onshore flooding is well modeled by SWE, which is fortunate since
the wetting-and-drying algorithms of GeoClaw can then be used.  In our
initial work we have simply suppressed the higher-order derivative terms
(switching to SWE) wherever the initial water depth was 10 m or less.
A better wave breaking model that allows dynamic switching has not yet been
implemented but will ultimately be incorporated.

Another issue to consider is the initial conditions for the simulation.
For the asteroid impact problem, even the Boussinesq equations are not
adequate to model the original generation or evolution of a deep crater
in the ocean.  We must start with the results of a three-dimensional
multiphysics hydrocode simulation and produce suitable initial conditions
for the depth-averaged equations.  We discuss this further in 
Section~\ref{sec:compResults}.

\section{Equations and Numerical Methods}\label{sec:math}

For simplicity we present the equations and algorithm primarily
in one space dimension and time,
since this is sufficient to illustrate the main ideas.

\subsection{The Shallow Water and Boussinesq Equations}

The shallow water equations can be written
\begin{eqnarray}\label{swec}
\begin{aligned}
h_t +(hu)_x &= 0\\
(hu)_t + \left(hu^2 + \half gh^2\right)_x &= -ghB_x,
\end{aligned}
\end{eqnarray}
where $h(x,t)$ = water depth, $B(x) =$ topography ~($B<0$ offshore),
$\eta(x,t)=B(x)+h(x,t)$, with  $\eta = 0$ being sea level.
Thus  $h_0(x) = -B(x)  $ is the depth of water at rest.
The depth-averaged horizontal velocity is $u(x,t)$, and so
$hu$ is the  momentum, and finally $g=9.81$ is the gravitional constant.
These equations are a long wave length approximation to the Euler
equations, in the limit of small ocean depth relative to the wavelength $L$ of the 
disturbance. For earthquake generated tsunamis, a typical ocean depth might be 4
km, and a subduction zone can have a wavelength of 100 km or more, giving a small ratio.

The equations \eqn{swec} have the form of a hyperbolic system of conservation
laws with a source term in the momentum equation that is nonzero only
on varying topography.
The GeoClaw implementation incorporates the topography term into the Riemann
solvers in order to obtain a well-balanced method \cite{LeVequeGeorgeEtAl2011},
which amounts to solving the equations \eqn{swec} in the nonconservative form
\begin{eqnarray}\label{swe}
\begin{aligned}
h_t +(hu)_x &= 0\\
(hu)_t + \left(hu^2 \right)_x + gh\eta_x&= 0.
\end{aligned}
\end{eqnarray}

Peregrine \cite{peregrine_long_1967} derived a Boussinesq-type extension
on a flat bottom in the form
\begin{eqnarray}
\begin{aligned}
h_t +(hu)_x &= 0\\
(hu)_t + \left(hu^2\right)_x + gh\eta_x
{ - \frac 1 3 h_0^2 (hu)_{txx}}&= 0
\end{aligned}
\end{eqnarray}
These equations have some drawbacks however, and do not match the dispersion
relation of the Euler equations as well as other models developed more 
recently. (For a historical review of Boussinesq-type models see 
\cite{brocchini:2013})

Madsen and Sorenson \cite{MadsenSorensen1992} and Shaffer and Madsen \cite{SchaefferMadsen1995}
optimized the equations by adding a term with a parameter $B_1$ that could
be chosen to match the
water wave dispersion relation more closely. On general topography,
they obtained
\begin{eqnarray}\label{SM}
\begin{aligned}
h_t +(hu)_x &= 0\\
(hu)_t + \left(hu^2\right)_x + gh\eta_x &=
\left(B_1 + \frac 1 2\right) h_0^2 (hu)_{txx} \\
&\qquad\null + \frac 1 6 h_0^3 (hu/h_0)_{txx}
- B_1 h_0^2 g(h_0\eta_x)_{xx} ,
\end{aligned}
\end{eqnarray}
where  $h_0(x)$ is the initial water depth, and matching the dispersion
relation leads to an optimal $B_1 = 1/15$. 

The equations \eqn{SM}
appear to have the form of the SWE \eqn{swe} together with source terms 
on the right hand side.  However, the standard fractional step approach
can't be used for equations in this form because the source term involves
$t$-derivatives.

These equations can be rewritten as
\begin{eqnarray} \label{bouss1}
\begin{aligned}
h_t +(hu)_x &= 0\\
(hu)_t + \left(hu^2\right)_x + gh\eta_x - { D_{11}((hu)_t) }
&= gB_1h_0^2(h_0\eta_x)_{xx}.
\end{aligned}
\end{eqnarray}
where the differential operator $D_{11}$ is defined by
\[
D_{11}(w) = (B_1 + 1/2)h_0^2 w_{xx} - \frac 1 6 h_0^3(w/h_0)_{xx}.
\]
Now subtracting 
$D_{11}\left(\left(hu^2\right)_x + gh\eta_x\right)$  from both sides 
of the momentum equation from \eqn{bouss1} gives
\begin{eqnarray}
\begin{aligned}
h_t +(hu)_x &= 0\\
{ [I - D_{11}]}\left[(hu)_t + \left(hu^2\right)_x + gh\eta_x\right] &=
{ -D_{11}\left[\left(hu^2\right)_x + gh\eta_x\right]} \\
&\qquad\null { + gh_0^2B_1(h_0\eta_x)_{xx}.}
\end{aligned}
\end{eqnarray}
By inverting $(I-D_{11})$ we get
\begin{eqnarray}\label{bouss2}
\begin{aligned}
h_t +(hu)_x &= 0\\
(hu)_t + \left(hu^2\right)_x + gh\eta_x &= \psi,
\end{aligned}
\end{eqnarray}
where $\psi$ is computed by solving an elliptic system:
\begin{equation}\label{elliptic}
[I - D_{11}] \psi = -D_{11}\left[\left(hu^2\right)_x + gh\eta_x\right]
+ gh_0^2B_1(h_0\eta_x)_{xx}.
\end{equation}
The system \eqn{bouss2} now looks like SWE plus a source term $\psi$ that
involves only spatial derivatives.

\subsection{Two dimensional versions}
For completeness, we include the two-dimensional version of these equations
to show that they have a similar structure.  In two dimensions let
$\vu = (u,v)$ be the two horizontal (depth-averaged) velocities. Then
the shallow water equations take the form
\begin{equation}\label{swe2v}
\begin{split} 
h_t + \grad\cdot (h\vu) &= 0\\
(h\vu)_t + \vz + gh \grad \eta &= 0.
\end{split} 
\end{equation}
where
\begin{equation}\label{z1}
\vz = \vu\,\grad\cdot(h\vu) + (h\vu\cdot\grad)\vu
= \bcm 
\left(hu^2\right)_x + (huv)_y \\
(huv)_x + \left(hv^2\right)_y \ecm.
\end{equation} 

The Boussinesq equations of \cite{MadsenSorensen1992,SchaefferMadsen1995}, 
as used in \cite{Kim2017}, take the form
\begin{equation}\label{bous2v}
\begin{split} 
h_t + \grad\cdot (h\vu) &= 0\\
(h\vu)_t + \vz + gh \grad \eta - D(h\vu)_t -gB_1 h_0^2 \grad\left(\grad
\cdot (h_0 \grad \eta)\right)&= 0.
\end{split} 
\end{equation}
where now the $2\times 2$ matrix $D$ consists of four linear 
differential operators:
\begin{equation}\label{D}
D = \bcm D_{11} & D_{12}\\
         D_{21} & D_{22} \ecm,
\end{equation}
with
\begin{equation}\label{Dij}
\begin{split} 
D_{11}(w) &= (B_1 + 1/2)h_0^2 w_{xx} - \frac 1 6 h_0^3(w/h_0)_{xx},\\
D_{12}(w) = D_{21}(w) &= (B_1 + 1/2)h_0^2 w_{xy} - \frac 1 6 h_0^3(w/h_0)_{xy},\\
D_{22}(w) &= (B_1 + 1/2)h_0^2 w_{yy} - \frac 1 6 h_0^3(w/h_0)_{yy}.
\end{split}
\end{equation}
As in 1D, in order to apply a fractional step method,
we subtract $D(\vz + gh \grad \eta)$ from both sides of the
momentum equation of \eqn{bous2v} so that it becomes
\begin{equation}\label{bous2mom}
[I-D][(h\vu)_t + \vz + gh \grad \eta] = -D(\vz + gh \grad \eta) 
+gB_1 h_0^2 \grad\left(\grad \cdot (h_0 \grad \eta)\right).
\end{equation}
Inverting $[I-D]$ allows rewriting \eqn{bous2v} as the SWE with a
source term,
\begin{equation}\label{bous2v2}
\begin{split} 
h_t + \grad\cdot (h\vu) &= 0\\
(h\vu)_t + \vz + gh \grad \eta &= \vpsi,
\end{split} 
\end{equation}
where $\vpsi$ involves only spatial derivatives and 
is determined by solving the elliptic equation
\begin{equation}\label{bous2psi}
[I-D]\vpsi  = -D(\vz + gh \grad \eta) 
+gB_1 h_0^2 \grad\left(\grad \cdot (h_0 \grad \eta)\right).
\end{equation}

Discretizing this elliptic equation leads
to a nonsymmetric sparse matrix with much
wider bandwidth than the tridiagonal matrix that arises in one dimension
due to the cross derivative terms.
Other than the significant increase in computing time required,
the fractional step and adaptive mesh refinement
algorithms described below carry over directly to the two dimensional
situation.

\subsection{Numerics}

We return to the one-dimensional equations in order to describe the numerical
algorithm and reformulation for patch-based adaptive refinement in space and
time.

We solve the one-dimensional Boussinesq equation \eqn{bouss2}, in which
$\psi$ is determined as the solution to \eqn{elliptic}, by
using a fractional step method with the following steps:
\begin{enumerate}
\item
Solve the elliptic equation \eqn{elliptic} for the source term $\psi$.
After this step the $\psi$ values are saved on each patch to use as 
boundary conditions for finer patches.
\item
Update the momentum by solving $(hu)_t = \psi$ over the time step
(e.g. with Forward Euler or two-stage Runge-Kutta method).
The depth $h$ does not change in this step since there is
no source term in the $h_t$ equation.
\item 
Take a step with the homogeneous SWE, using the results of step 2 as initial
data.  This step uses the regular GeoClaw software and Riemann solvers.
\end{enumerate}

We solve the implicit system first and then take the shallow
water step because this facilitates interpolating in time for values required
on the edge of grid patches.
In order to explain this in more detail, we
introduce some notation for a simple case in one space dimension.

First suppose we only have a single grid at one resolution, with no AMR.
We denote the numerical solution at some time $t_N$ by
$(H,HU)^N$, the cell-averaged approximations to
depth and momentum on the grid.  We also use 
$\Psi^N$ for the source
term determined by solving the discrete elliptic system defined by $(H,HU)^N$
on this grid. 
We also assume at the start of the time step that we have boundary conditions for
$(H,HU)$ and also for the Boussinesq correction $\Psi$, provided in the form
of ``ghost cell'' values in a layer of cells surrounding the grid (or two
layers in the case of $(H,HU)$ since the high-resolution explicit methods 
used have a stencil of width 5 because of slope limiters).
On a single
grid we assume that it is sufficient to use the Dirichlet condition $\Psi=0$
in all ghost cells surrounding the grid, i.e. that there is no Boussinesq
correction in these cells. This is reasonable for a large domain where the
waves of interest are confined to the interior of the domain.  We also use
zero-order extrapolation boundary condition for $(H,HU)$, which give a
reasonable non-reflecting boundary condition for the 
SWE step \cite{LeVequeGeorgeEtAl2011}.

A single time step of the fractional step algorithm on this grid then takes
the following form in order to advance 
$(H,HU)^N$ to $(H,HU)^{N+1}$ at time $t_{N+1} = t_N + \dt$:
\begin{enumerate}
\item Solve the elliptic system for $\Psi^N$.  The right-hand side depends on
$(H,HU)^N$.
\item Advance the solution using the source terms (Boussinesq corrections):\\
$H^* = H^N,~ (HU)^* = (HU)^N + \dt \Psi^N$ (using Forward Euler, for example).
\item Take a time step of length $\dt$ with the SWE solver, with initial data
$(H,HU)^*$, to obtain $(H,HU)^{N+1}$.
We denote this by $(H,HU)^{N+1} = SW((H,HU)^*, \dt)$.
\end{enumerate}
In the software it is convenient to store the source term at each time as
another component of the solution vector, so we also use $Q^N=(H,HU,\Psi)^N$
to denote this full solution at time $t_N$.

Now suppose we
have two grid levels with refinement by a factor of 2 in time. 
We denote the coarse grid values at some time $t_N$ as above.
We assume that the fine grid is at time $t_N$, but that
on the fine grid we must take two time steps of $\dt/2$ to reach time
$t_{N+1}$.  We denote the fine grid values at time $t_N$ using
lower case, $(h,hu)^N$ and $q^N=(h,hu,\psi)^N$.
We also need boundary conditions in the ghost cells of the fine grid patch.
If a patch edge is coincident with a domain boundary, then we use the 
Dirichlet BC $\psi=0$ and extrapolation BCs for $(h,hu)$, as described above.
For ghost cells that are interior to the coarse grid, we let 
$\If(Q)$ represent a spatial interpolation operator that 
interpolates from coarse grid values to the ghost cells of a fine grid
patch at time $t_N$.  This operator is applied to all three components of
$Q^N$, i.e. to the source term as well as the depth and momentum, in order
to obtain the necessary boundary conditions for $q^N$.

Then one time step on the coarse grid, coupled with two time steps on the 
fine grid, is accomplished by the following steps:


\begin{enumerate}

\item Coarse grid step:
    \begin{enumerate}
    \item Take time step $\Delta t$ on the coarse grid as described above for
    the single grid algorithm, but denote the result by 
    $(\widetilde H,\widetilde{HU})^{N+1}$ since these provisional values will later
    be updated.
    \item Using the Dirichlet BCs $\Psi=0$ on the domain boundary, solve for
    a provisional $\widetilde\Psi^{N+1}$.  This will be needed for interpolation
    in time when determining boundary conditions for $\psi$ on the fine grid,
    using $\If(Q^N)$ and $\If(\widetilde Q^{N+1})$.
    \end{enumerate}
    
\item Fine grid steps:
    \begin{enumerate}
    \item Given $(h,hu)^N$ and boundary conditions $\If(Q^N)$, solve the
    elliptic system for $\psi^N$.
    \item Update using the source terms, $(h,hu)^* = (h,hu)^N + (0, \frac{\dt}{2}\psi^N)$.
    \item Take a shallow water step: $(h,hu)^{N+1/2} = SW((h,hu)^*, \dt/2)$.\\
    Note that we use $t_{N+1/2} = t_N + \dt/2$ to denote the intermediate time.
    \item Obtain BCs at this intermediate time as 
    $\half(\If(Q^N)+\If(\widetilde Q^{N+1}))$.
    \item Solve the elliptic system for $\psi^{N+1/2}$.
    \item Update using the source terms, $(h,hu)^* = (h,hu)^{N+1/2} + (0, \frac{\dt}{2}\psi^{N+1/2})$.
    \item Take a shallow water step: $(h,hu)^{N+1} = SW((h,hu)^{*}, \dt/2)$.
    \end{enumerate}
    
\item Update coarse grid:
    \begin{enumerate}
    \item Define $(H,HU)^{N+1}$ by the provisional values by $(\widetilde H,\widetilde{HU})^{N+1}$ where there is no fine grid covering a grid cell,
    but replacing
    $(\widetilde H,\widetilde{HU})^{N+1}$ by the average of $(h,hu)^{N+1}$ 
    over fine grid cells that cover any coarse grid cell.
    \end{enumerate}
\end{enumerate}

The final step is applied because the fine grid values $(h,hu)^{N+1}$ are
more accurate than the provisional coarse grid values.

We then proceed to the next coarse grid time step.  
Note that at the start of this step, the updated $(H,HU)^{N+1}$ will be used
to solve for $\Psi^{N+1}$.  The provisional $\widetilde \Psi^{N+1}$ is discarded.
Hence two elliptic solves are required on the coarse level each time step,
rather than only one as in the single grid algorithm.

If the refinement factor is larger than 2, then the same approach outlined
above works, but there will be additional time steps on level 2.  For each
time step the ghost cell BCs will be determined by linear interpolation in 
time between $\If(Q^N)$ and $\If(\widetilde Q^{N+1})$.

If there are more than two levels, then this same idea is applied recursively:
After each time step on level 2, any level 3 grids will be advanced by the
necessary number of time steps to reach the advanced time on level 2.  In this
case there will also be two elliptic solves for every time step on level 2,
once for the provisional values after advancing level 2, and once at the start
of the next level 2 time step after $(h,hu)$ on level 2 has been updated by
averaging the more accurate level 3 values.

\section{Computational Results}\label{sec:compResults}
In this section we show an end-to-end simulation of a hypothetical asteroid impact off
the coast of Washington, from initial conditions to shoreline inundation. 
We present our initialization procedure in some detail. 
The section ends with a discussion of results.

\subsection{Initialization Procedure}\label{subsec:init}
\begin{figure}[t]
(a) \includegraphics[height=1.8in]{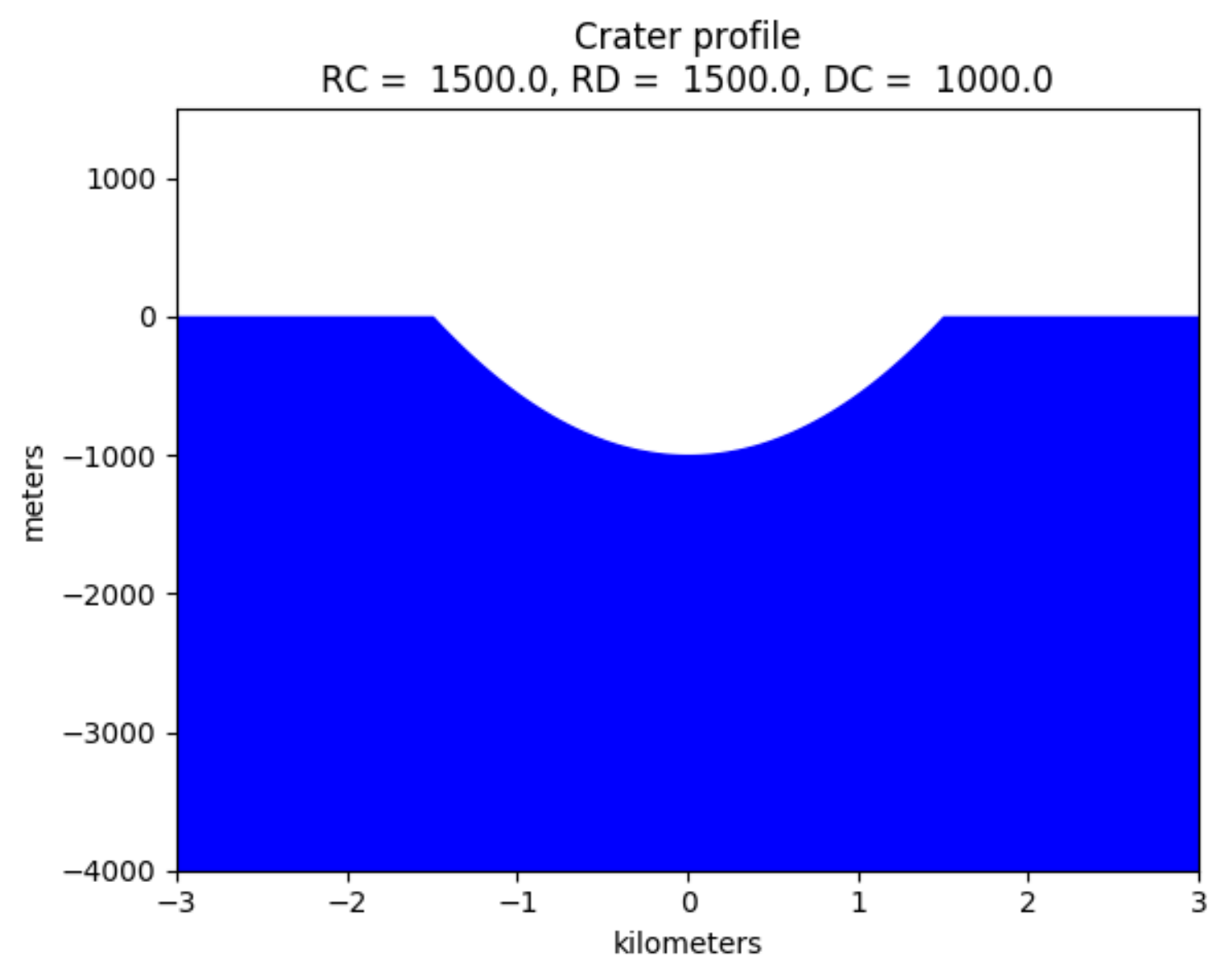}
\hspace*{.3in}
(b) \includegraphics[height=1.8in]{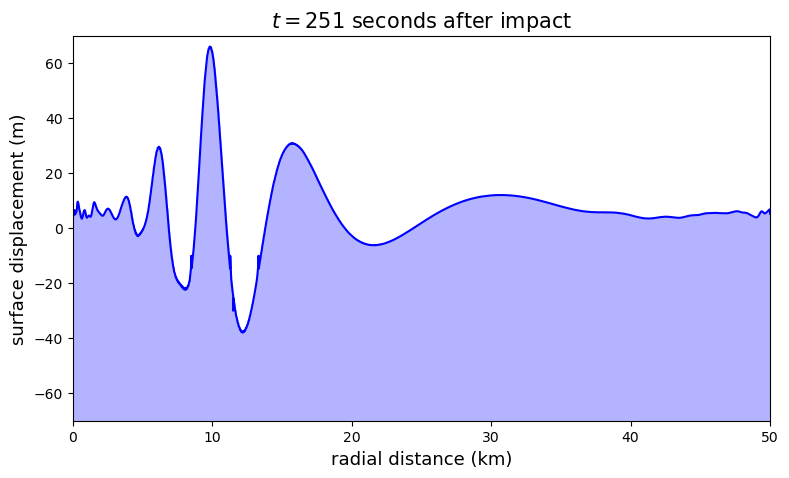}
\caption{(a) initial conditions of a static crater of depth 1
km.~and diameter 3 km., which were the initial conditions for the hydrocode
simulation. (b) The radially symmetric results of the hydrocode
simulation at 251 seconds, used to start the GeoClaw Boussinesq
simulation. \label{fig:craterfig}}
\end{figure}

The computational results presented in this section use initial conditions of a 
static crater, illustrated in Figure \ref{fig:craterfig}(a). This is a standard
test problem in the literature, from \cite{WardAsphaug:2000}.  The crater is one 
kilometer deep, with a diameter of 3 kilometers, 
in an ocean of depth 4 kilometers. Depth-averaged equations are unsuitable 
for modeling the generation of the crater and the initial flow, 
since the ensuing large vertical velocity 
components are not modeled. The initial conditions for our depth-averaged
simulations are taken from a three-dimensional 
hydrocode simulation of the first 251 seconds after impact. The 
hydrocode ALE3D \cite{ALE3d} was run by a collaborator \cite{RobertsonPC}
in a radially symmetric manner, and the surface 
displacement at time $t=251$ seconds was recorded, as shown in Figure
\ref{fig:craterfig}(b). It proved
too noisy to depth-average the horizontal velocity from the hydrocode. 
Instead we set the velocity based on the surface displacement and assuming
that the wave was a purely outgoing wave satisfying the SWE, for which the
velocity then depends only on the ocean depth and surface displacement.
We could then place the ``initial'' crater anywhere in the ocean.

This procedure for initialization of the velocity can be done because
the wave speed for the SWE is independent of wave number, and the eigenvectors
of the linearized Jacobian matrix give the relation between surface elevation
and fluid velocity for unidirectional waves.  However, in the
Boussinesq equations the wave speed depends on wave number and using the
initialization based on the SWE results in a small wave propagating inward as
well.  Better initialization procedures will be investigated in future research.

\subsection{Adaptive simulation}\label{subsec:adapt}

\begin{figure}[t]
\centerline{\includegraphics[width=0.75\textwidth]{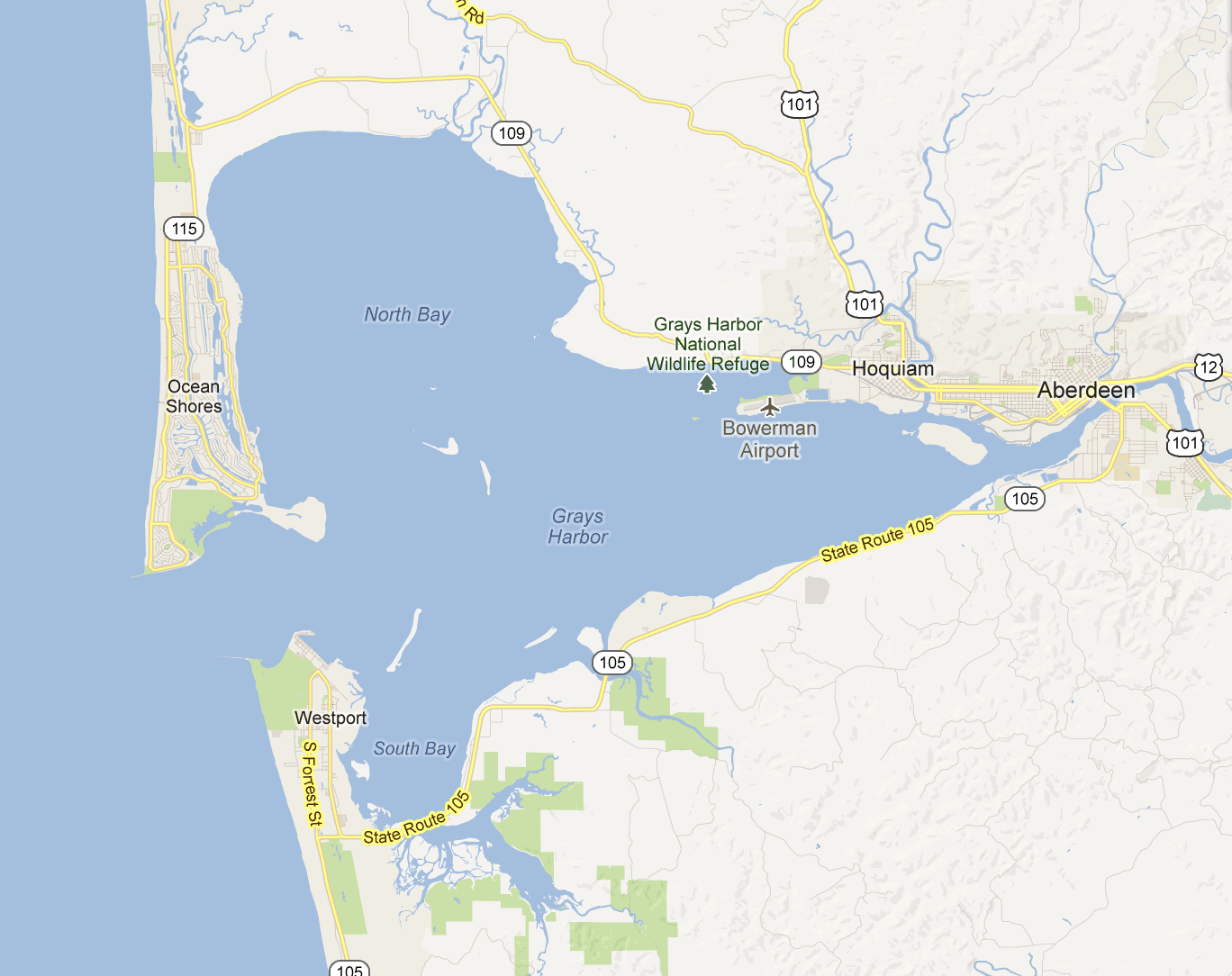}}
\caption{Figure shows a Google Maps screenshot of
Grays Harbor, on the Washington coast. The
community of Westport is on the southern peninsula.
This is the focus of the inundation modeling
presented in Figure~\ref{GHFigs}.
\label{fig:grayFig}}
\end{figure}

As an illustration, we show a simulation of a hypothetical asteroid impact off the
coast of Washington. We place the initial crater approximately 150 km west of 
Grays Harbor, to study the vulnerable area around Westport, WA, shown in
Figure \ref{fig:grayFig}. This well-studied area is in close proximity to the 
Cascadia Subduction Zone, which can generate Mw 9 earthquakes.  
The Ocosta elementary school in Westport was recently rebuilt to 
incorporate the first tsunami vertical evacuation structure in the U.S., 
due to the low topography of this region, with design work based in part on
GeoClaw modeling \cite{Ocosta}.  Detailed bathymetry and topography data is
available in this region at a resolution of 1/3 arcsecond \cite{Astoria-DEM},
which is roughly 10 meters in latitude and 7 m in longitude at this location.
For the ocean we use the etopo1 topography DEM \cite{etopo1} at a resolution
of 1 arcminute.

A 7 level simulation was used, starting with a coarsest level over the ocean
with $\Delta y = 10$ arcminutes, and refining by factors 5, 3, 2, 2, 5, 6
at successive levels, with an overall refinement by 1800 for the level-7 grids
with $\Delta y = 1/3$ arcsecond.  On each grid $\Delta x = 1.5\Delta y$ so
that the finest-level computational grids are at a resolution of
roughly 10 m in both $x$ and $y$. 

\begin{figure}[th]

\includegraphics[width=0.45\textwidth]{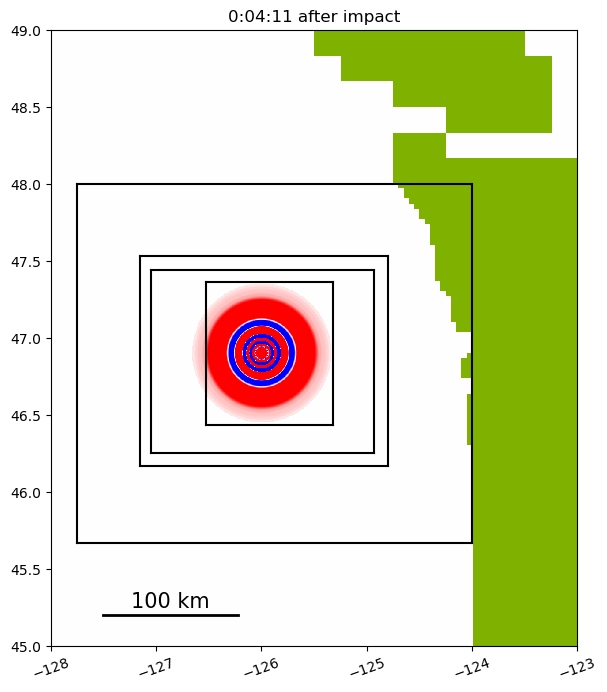}
\hspace{5pt}
\includegraphics[width=0.45\textwidth]{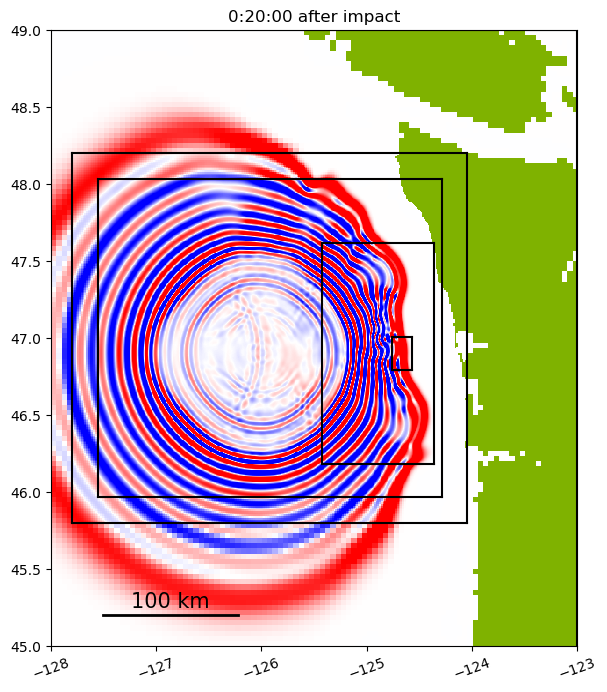}

\includegraphics[width=0.45\textwidth]{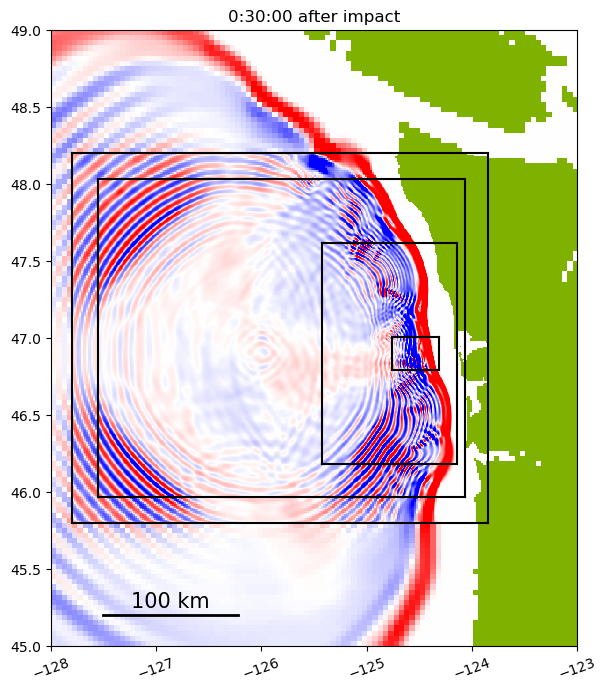}
\hspace{5pt}
\includegraphics[width=0.45\textwidth]{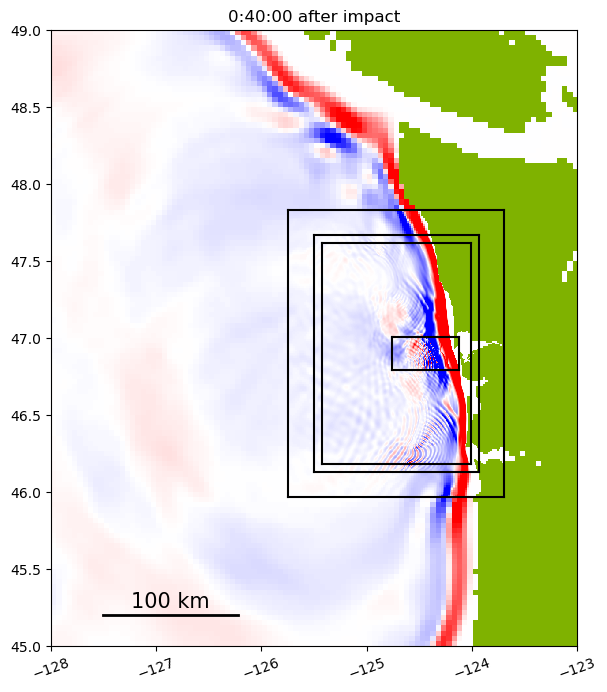}
\caption{Figure shows initial conditions and mesh configuration and
three later times during the adaptive solution.  The black rectangles show
boundaries of refined patches with finer resolution.  The colors show elevation
above (red) or below (blue) sea level, and saturate at $\pm 3$ meters.
The waves are larger amplitude but this color range is used to also show
the smaller waves in the oscillatory wave train.  Note that these are
not resolved on coarser levels, and that AMR is guided to focus on the
waves approach Grays Harbor, WA.
\label{firstSetCompFigs}}
\end{figure}

\begin{figure}[th]

\includegraphics[width=0.95\textwidth]{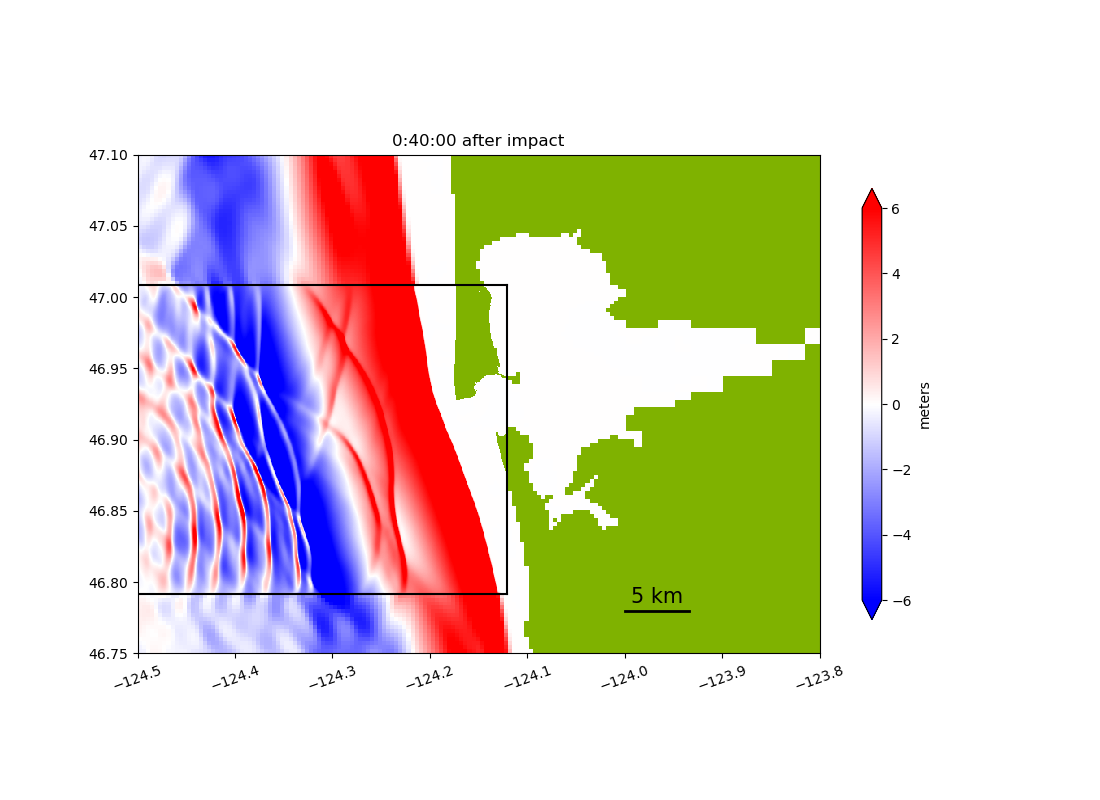}

\includegraphics[width=0.95\textwidth]{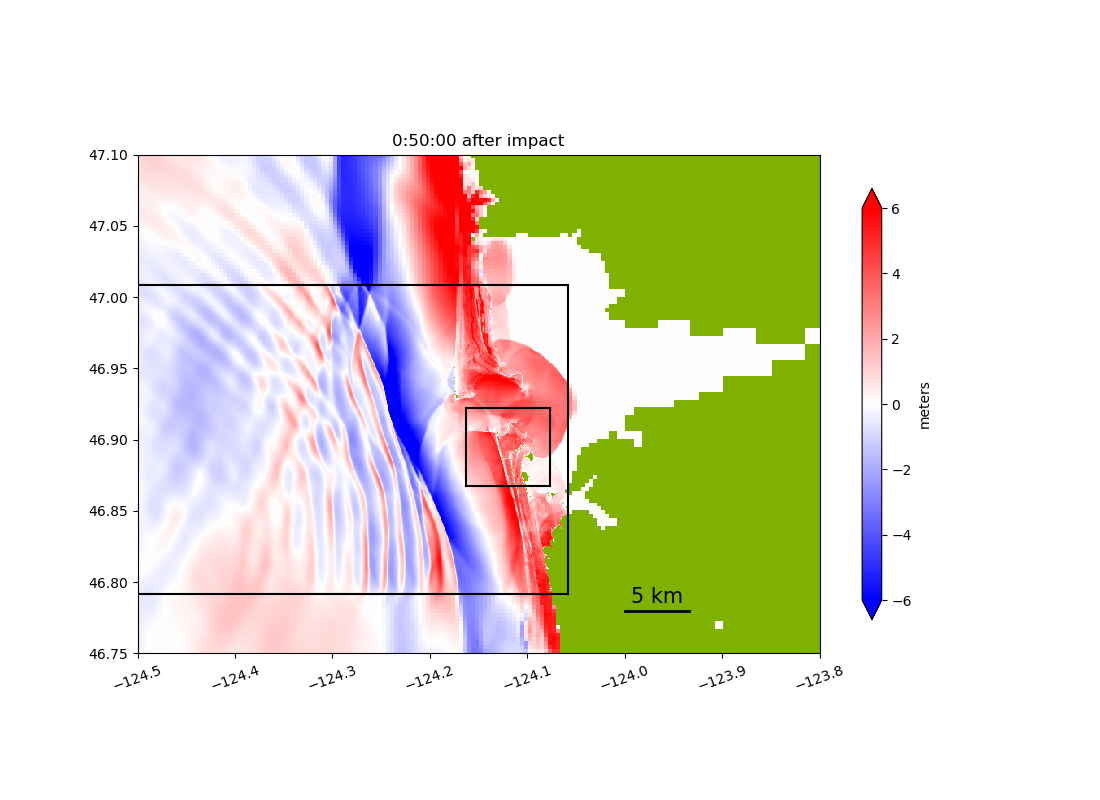}

\caption{The tsunami approaching Grays Harbor and overtopping the
Ocean Shores and Westport peninsulas.  
In this figure, the colors saturate at $\pm 6$ meters elevation relative
to sea level.
\label{GHFigs}}
\end{figure}

Figures \ref{firstSetCompFigs} and \ref{GHFigs} 
show snapshots of the simulation at the
indicated times. The dispersion is clearly evident, with a much more
oscillatory solution than would be obtained with the shallow water
equations.
The figures show how the fine grid patches move to follow the 
expanding wave. The refinement is guided to focus at later times only on waves
approaching Grays Harbor.  There are some 
reflections at the grid boundaries, but they are much smaller in magnitude 
than the waves we are tracking. The  6th level refined patch appears 
approximately 20 minutes after the impact, to track the waves as they 
approach Grays Harbor. Note how the bathymetry is refined along with the 
solution when the finer patches appear.  The close up plots near Grays Harbor
in Figure~\ref{GHFigs}
shows ``soliton fission'', a nonlinear dispersive wave phenomenon seen 
near the coast that can be captured with Boussinesq solvers
\cite{BabaTakahashiEtAl2015,MatsuyamaIkenoEtAl2007}.
In this simulation, we switch from Boussinesq to SWE at a depth of 10 meters
(based on the undisturbed water depth).
Figure~\ref{GHFigs} also shows the waves sweeping over the Westport and
Ocean Shores peninsulas. In this calculation, the finest level 7 grid was
placed only over Westport, while the level 6 grid covers both peninsulas.

\subsection{Discussion of Results}\label{subsec:discuss}

Since the wavelengths of asteroid-generated tsunamis are shorter than 
those of earthquake-generated tsunamis, dispersive effects may be very
important. Dispersion gives rise to a highly oscillatory set of waves that
propagate at different speeds, possibly affecting the arrival time of the
first wave and leading to significant waves over a longer time period.
These waves can undergo substantial amplification during the shoaling
process on the continental shelf and break up into large solitary waves.

On the other hand, wave breaking can rapidly dissipate the energy in
short wavelength waves.  Moreover, a train of waves approaching shore 
results in nonlinear interactions in the swash zone, where waves run up the
beach. In the swash zone, a large approaching
wave may be largely negated by the rundown of the previous wave.  Thus it
is not clear {\em a priori} whether waves modeled with the Boussinesq equations
will result in substantially different onshore inundation than would be
observed is only using the SWE, for which computations are much less expensive.  
This question still remains to be answered.

\clearpage

\section{Open Problems and Future Research}\label{sec:conc}
We have demonstrated that it is possible to develop a high-fidelity 
modeling capability that includes propagation and inundation by combining the 
Boussinesq equations, shallow water equations, and patch-based adaptive 
mesh refinement in space and time.  We are embedding this in the GeoClaw software 
framework, which has previously been well-validated for earthquake-generated 
tsunamis.  This will allow the efficient simulation of dispersive waves
generated from asteroid impacts as they propagate across the ocean, combined
with high-resolution simulation of the resulting inundation on the coast.
This capability could be very important in hazard assessment for an incipient
impact.

Unfortunately the software is not yet robust enough for general use. While it
often works well, stability issues sometimes arise at the edges of finer 
grid patches when refinement ratios greater than 2 are used
from one level to the next. Larger refinement ratios are generally desirable,
so that an overall refinement factor of several thousand between the 
coarsest and finest meshes can be achieved with 6 or 7  levels of mesh refinement.  
Other Boussinesq codes we are aware of use
factor of 2 refinement \cite{Popinet:2015,ForestClaw:2017}, so this may be an 
inherent instability.  Moreover,
numerical instabilities have also been observed by other researchers
\cite{LovholtPedersen2009}
even on a uniform grid when there are sharp changes in bathymetry rather
than in the grid resolution.  We are investigating this issue
theoretically, and may find that a different discretization or even a 
different fomulation of the Boussinesq equations is required to obtain a 
sufficiently robust code.

We are also continuing to investigate the shoaling phenomena and the best
way to incorporate wave breaking and the transition from Boussinesq to shallow
water equations.   This can have a significant impact
on the resulting onshore runup and inundation.

Not discussed here but currently under investigation is the possibility of 
solving the implicit system of  equations on multiple levels in a 
coupled manner, whenever
the levels have been advanced to the same point in time.  It is an
open question whether a coupled system is more accurate and/or stable, 
and possibly less computationally expensive, in the context of patch-based
adaptive mesh refinement that includes refinement in time.

Finally, as mentioned in section \ref{subsec:init}, additional research is 
needed on ways to initialize the depth-averaged model.  
Better initialization procedures will allow a more seamless transition
from three-dimentional hydrocode simulations of asteroid impacts in the ocean
to our model of tsunami propagation and inundation.

\vskip 5pt

{\bf Acknowledgments.}
We thank Darrel Robertson, a member of the ATAP team, for providing the
hydrocode simulation results for our initial conditions.

This work was partially supported by the NASA Asteroid Threat Assessment
Project (ATAP) through the Planetary Defense Coordination Office and
BAERI contract AO9667.


\begin{thebibliography}{10}

\bibitem{etopo1}
C.~Amante and B.~W. Eakins.
\newblock {ETOPO1 1 Arc-Minute Global Relief Model: Procedures, Data Sources
  and Analysis}.
\newblock NOAA Technical Memorandum NESDIS NGDC-24,
  http://www.ngdc.noaa.gov/mgg/global/global.html, 2009.

\bibitem{BabaTakahashiEtAl2015}
Toshitaka Baba, Narumi Takahashi, Yoshiyuki Kaneda, Kazuto Ando, Daisuke
  Matsuoka, and Toshihiro Kato.
\newblock Parallel implementation of dispersive tsunami wave modeling with a
  nesting algorithm for the 2011 {Tohoku} tsunami.
\newblock {\em Pure and Applied Geophysics}, 172(12):3455--3472, 2015.

\bibitem{BergerGeorgeEtAl2011a}
M.~J. Berger, D.~L. George, R.~J. {LeVeque}, and K.~T. Mandli.
\newblock The {GeoClaw} software for depth-averaged flows with adaptive
  refinement.
\newblock {\em Advances in Water Resources}, 34(9):1195--1206, September 2011.

\bibitem{NASA-AGT-2018}
M.~J. Berger and R.~J. LeVeque.
\newblock Modeling issues in asteroid-generated tsunamis.
\newblock NASA Technical Memorandum NASA/CR-2018-219786, ARC-E-DAA-TN53167,
  http://hdl.handle.net/2060/20180006617, 2018.

\bibitem{pardiso-7.2b}
Matthias Bollh{\"o}fer, Olaf Schenk, Radim Janalik, Steve Hamm, and Kiran
  Gullapalli.
\newblock State-of-the-art sparse direct solvers.
\newblock In Ananth Grama and Ahmed~H. Sameh, editors, {\em Parallel Algorithms
  in Computational Science and Engineering}, pages 3--33, Cham, 2020. Springer
  International Publishing.

\bibitem{brocchini:2013}
M.~Brocchini.
\newblock A reasoned overview on {Boussinesq-type} models: the interplay
  between physics, mathematics and numerics.
\newblock {\em Phil. Trans. Royal Soc. A}, 469(20130496), 2013.

\bibitem{ForestClaw:2017}
D.~Calhoun and C.~Burstedde.
\newblock Forestclaw: A parallel algorithm for patch-based adaptive mesh
  refinement on a forest of quadtrees, 2017.
\newblock ArXiv:1703.03116.

\bibitem{clawpack}
{Clawpack Development Team}.
\newblock Clawpack software.
\newblock http://www.clawpack.org, 2021.

\bibitem{Glimsdal2013}
S.~Glimsdal, G.~K. Pedersen, C.~B. Harbitz, and F.~L{\o}vholt.
\newblock Dispersion of tsunamis: Does it really matter?
\newblock {\em Nat. Hazards Earth Syst. Sci.}, 13:1507--1526, 2013.

\bibitem{geoclaw-nthmp-results:2011}
F.~{Gonz\'alez}, R.~J. LeVeque, J.~Varkovitzky, P.~Chamberlain, B.~Hirai, and
  D.~L. George.
\newblock {GeoClaw Results for the NTHMP Tsunami Benchmark Problems}.
\newblock http://depts.washington.edu/clawpack/links/nthmp-benchmarks, 2011.

\bibitem{Ocosta}
F.~I. Gonzal\'ez, R.~J. LeVeque, and L.~M. Adams.
\newblock {Tsunami Hazard Assessment of the Ocosta School Site in Westport,
  WA}.
\newblock http://hdl.handle.net/1773/24054, 2013.

\bibitem{Kim2017}
J.~Kim, G.~K. Pedersen, F.~L{\o}vholt, and R.~J. LeVeque.
\newblock A {Boussinesq} type extension of the {GeoClaw} model - a study of
  wave breaking phenomena applying dispersive long wave models.
\newblock {\em Coastal Engineering}, 122:75 -- 86, 2017.

\bibitem{KorycanskyLynett}
D.~G. Korycansky and P.~J. Lynett.
\newblock Offshore breaking of impact tsunami: {{The Van Dorn}} effect
  revisited.
\newblock {\em Geophysical Research Letters}, 32(10), 2005.

\bibitem{KorycanskyLynett2007}
D.~G. Korycansky and Patrick~J. Lynett.
\newblock Run-up from impact tsunami.
\newblock {\em Geophysical Journal International}, 170:1076--1088, 2007.

\bibitem{fvmhp}
R.~J. LeVeque.
\newblock {\em Finite Volume Methods for Hyperbolic Problems}.
\newblock Cambridge University Press, 2002.

\bibitem{LeVequeGeorgeEtAl2011}
R.~J. LeVeque, D.~L. George, and M.~J. Berger.
\newblock Tsunami modelling with adaptively refined finite volume methods.
\newblock {\em Acta Numerica}, 20:211 -- 289, May 2011.

\bibitem{LovholtPedersen2009}
F.~L{\o}vholt and G.~Pedersen.
\newblock {Instabilities of Boussinesq models in non-uniform depth}.
\newblock {\em International Journal for Numerical Methods in Fluids},
  61:606--637, 2009.

\bibitem{MadsenSorensen1992}
P.~A. Madsen and O.~R. S{\o}rensen.
\newblock A new form of the {{Boussinesq}} equations with improved linear
  dispersion characteristics. {{Part}} 2. {{A}} slowly-varying bathymetry.
\newblock {\em Coastal engineering}, 18(3-4):183--204, 1992.

\bibitem{MatsuyamaIkenoEtAl2007}
Masafumi Matsuyama, Masaaki Ikeno, Tsutomu Sakakiyama, and Tomoyoshi Takeda.
\newblock A study of tsunami wave fission in an undistorted experiment.
\newblock {\em Pure and Applied Geophysics}, 164(2):617--631, 2007.

\bibitem{ALE3d}
A.~L. Nichols and D.~M.~Dawson (Eds.).
\newblock {ALE3D User’s Manual: An Arbitrary Lagrange/Eulerian 2D and 3D Code
  System}.
\newblock Lawrence Livermore National Laboratory LLNL-SM-726137, 2017.

\bibitem{Astoria-DEM}
{NOAA National Geophysical Data Center}.
\newblock {Astoria, Oregon 1/3 Arc-second MHW Coastal Digital Elevation Model}.
\newblock
  https://www.ncei.noaa.gov/access/metadata/landing-page/bin/iso?id=gov.noaa.ngdc.mgg.dem:5490,
  Accessed 2019.

\bibitem{peregrine_long_1967}
D.~H. Peregrine.
\newblock Long waves on a beach.
\newblock {\em Journal of Fluid Mechanics}, 27(4):815--827, 1967.

\bibitem{Popinet:2015}
S.~Popinet.
\newblock A quadtree-adaptive multigrid solver for the {S}erre-{G}reen-{N}aghdi
  equations.
\newblock {\em J. Comp. Phys.}, 302:336--358, 2015.

\bibitem{RobertsonPC}
D.~Robertson.
\newblock {Personal Communication}, 2019.

\bibitem{RobertsonGisler2019}
Darrel~K. Robertson and Galen~R. Gisler.
\newblock Near and far-field hazards of asteroid impacts in oceans.
\newblock {\em Acta Astronautica}, 156:262--277, 2019.

\bibitem{SchaefferMadsen1995}
H.~A. Sch{\"a}ffer and P.~A. Madsen.
\newblock Further enhancements of {{Boussinesq}}-type equations.
\newblock {\em Coastal Engineering}, 26(1-2):1--14, 1995.

\bibitem{WardAsphaug:2000}
S.~N. Ward and E.~Asphaug.
\newblock Asteroid impact tsunami: A probabilistic hazard assessment.
\newblock {\em Icarus}, 145:64--78, 2000.

\end{thebibliography}

\end{document}